\documentclass{article}
\usepackage{epsfig,amsmath,amsfonts}

\newcommand{\C}{\mathbb C}  
  
\newcommand{\PP}{\mathbb P} 
\newcommand{\R}{\mathbb R}

\newcommand{\G}{\mathbb G}
\newcommand{\V}{\mathbb V}

\newcommand{\rank}{\mbox{rank}}
\newcommand{\grad}{\mbox{grad}}
\newcommand{\Perm}{\mbox{Perm}}
\newcommand{\im}{\mbox{im}}

\newtheorem{lemma}{Lemma}
\newtheorem{theorem}{Theorem}

\newtheorem{definition}{Definition}
\newtheorem{prop}{Proposition}

\newcommand{\join}{\vee}
\newcommand{\meet}{\wedge}

\newcommand{\ra}{\rightarrow}

\newcommand{\lra}{\longrightarrow}

\newcommand{\al}{\alpha}

\newcommand{\ga}{\gamma}

\newcommand{\ep}{\epsilon}
\newcommand{\la}{\lambda}
\newcommand{\om}{\omega}
\newcommand{\si}{\sigma}
\newcommand{\tht}{\theta}

\newcommand{\Ga}{\Gamma}
\newcommand{\De}{\Delta}

\newcommand{\La}{\Lambda}
\newcommand{\Om}{\Omega}

\newcommand{\Th}{\Theta}

\newcommand{\wtilde}{\widetilde}
\newcommand{\what}{\widehat}

\newcommand{\nodes}{\mbox{nodes}}

\newenvironment{proof}{{\em Proof:}}{\hfill\rule{2mm}{2mm}}
\begin{document}

\title{\bf {Recovering an Algebraic Curve Using its Projections From
Different Points} 
{Applications to Static and Dynamic Computational Vision}}
\author{Michael Fryers$^1$, Jeremy Yirmeyahu Kaminski$^2$ and Mina Teicher$^2$
\thanks{This work is partially supported by the  Emmy Noether Institute for
Mathematics and the Minerva Foundation of Germany, by the Excellency Center
of the Israel Science Foundation "Group Theoretic Methods in the Study of
Algebraic Varieties" and by EAGER (European network
in Algebraic Geometry).}\\
$^1$ Universit\"{a}t Hannover, Hannover, Germany. e-mail: \\
fryers@math.uni-hannover.de \\
$^2$ Bar-Ilan University, Ramat-Gan, \\
Israel. e-mail:\{kaminsj,teicher\}@macs.biu.ac.il }

%\date{}

\maketitle

\begin{abstract}
We study how an irreducible closed algebraic curve $X$ embedded in
$\C \PP^3$, which degree is $d$ and genus $g$, can be recovered using its
projections from points onto embedded projective planes. The different
embeddings are unknown. The only input is the defining equation of
each projected curve. We show how both the embeddings and the curve in
$\C \PP^3$ can be recovered modulo some actions of the group of projective
transformations of $\C \PP^3$.  

In particular in the case of two projections, we show how in a generic
situation, a characteristic matrix of the two embeddings can be recovered. In the
process we address dimensional issues and as a result establish the
minimal number of irreducible algebraic curves required to compute
this characteristic matrix up to a finite-fold ambiguity, as a function of their
degree and genus. Then we use this matrix to recover the class of the
couple of maps and as a consequence to recover the curve. For
a generic situation, two projections define a curve with two irreducible
components. One component has degree $d(d-1)$ and the other has degree
$d$, being the original curve.  

Then we consider another problem. $N$ projections, with known
projections operators and $N >> 1$, are considered as an input and we
want to recover the curve. The recovery can be done by linear
computations in the dual space and in the Grassmannian of lines in
$\C \PP^3$, that we denote by $\G(1,3)$. Those computations are
respectively based on the dual variety and on the variety of intersecting lines. In both
cases a simple lower bound for
the number of necessary projections is given as a function of degree
and genus. A closely related question is also considered. Each point
of a finite closed subset of an irreducible algebraic curve, is
projected onto a plane from a point. For each point the center of
projection is different. The projections operators are known. We show
when and how the recovery of the algebraic curve is possible in 
function of the degree of the curve of minimal degree
generated by the centers of projection.

These questions were motivated by applications to static and dynamic
computational vision. Therefore a second part of this work is devoted
to applications to this field. The results in this paper solve a long
standing problem in computer vision that could not have been solved
without algebraic-geometric methods. 
\end{abstract}

%=====================
%=====================
\section{Introduction}
\label{sec::intro}

Consider an irreducible closed algebraic curve $X \in \C\PP^3$ (in the
sequel we simply write $\PP^n$ for $\C\PP^n$). This curve is projected
onto several projective planes embedded in $\PP^3$ through several
center of projections, say $\{{\bf O}\}_i, i=1,..,n$. Each projection
mapping, denoted by $\pi_i:\PP^3 \setminus \{{\bf O}_i\} \lra \PP^2$ is
presented as a $3 \times 4$ matrix, ${\bf M}_i$ defined modulo
multiplication by non-zero scalar. Then each point ${\bf
P}$ different from ${\bf O}_i$ is mapped by $\pi_i$ to 
${\bf M}_i {\bf P}$. Each projection operator $\pi_i$, via its matrix,
can be regarded as a point in
$\PP^{11}$. Let $Y_i = \pi_i(X)$ be the different 
projections of the curve $X$. In the sequel we always deal with
generic configurations, even when not mentionned explicitly.

When we consider the problem of recovering the projections maps from
the projected curves, we will show that the recovery is possible only
modulo some action of the group of projective transformations of $\PP^3$
on the set of projection maps. To define 
this action we refer to a projection map as a point in
$\PP^{11}$. Assume that we have $n$ projections, consider the
following projective variety $\V = \stackrel{n \mbox{ } 
times}{\overbrace{\PP^{11} \times ... \times \PP^{11}}}$. Let $Pr_3$ be the group of
projective transformations of $\PP^3$. We define an action of $\Pr_3$
on $\V$ as follows: $\tht_n: Pr_3 \lra Mor(\V,\V), {\bf A} \mapsto (({\bf
Q}_1,...,{\bf Q}_n) \mapsto ({\bf M}_1 {\bf A}^{-1},...,{\bf
M}_n {\bf A}^{-1}))$, where each matrix ${\bf M}_i$ is built from
the coordinates of ${\bf Q}_i = [Q_{i,1},...,Q_{i,12}]^T$ as follows:
$$
{\bf M}_i = \left[\begin{array}{cccc}
		Q_{i,1} & Q_{i,2} & Q_{i,3} & Q_{i,4} \\
		Q_{i,5} & Q_{i,6} & Q_{i,7} & Q_{i,8} \\
		Q_{i,9} & Q_{i,10} & Q_{i,11} & Q_{i,12}
	    \end{array}\right]
$$ 
The geometric meaning of this action is that if we change the
projective basis in $\PP^3$, by the transformation ${\bf A}$, we need to
change the projection maps in accordance for the projected curves to
be invariant.

We first investigate the case of two projections. Given the projected
curves $Y_1$ and $Y_2$ as the only data, our first problem is to
compute the characteristic matrix (to be defined below) of two
projections maps, $\pi_1$ and $\pi_2$, up to a 
finite-fold ambiguity. It is shown that this is equivalent to
find a necessary and sufficient conditions on $X$ for the action of
$\tht_2$ to have a finite number of orbits. Then we show that for each orbit we
can recover the curve $X$ modulo $Pr_3$. More precisely each orbit
induces a curve embedded in $\PP^3$ containing two irreducible
components, one of degree $d(d-1)$ and the other of degree $d$. The
latter is the curve we are looking for. 

Then we turn to another problem. The projections maps $\pi_i,
i=1,..,n$ are now assumed to be known, in addition to the projected
curves $Y_i$. We want to recover the curve $X$. This can be performed
by linear computations using either the dual varieties or the variety
of lines intersecting $X$. In both case a simple lower bound of the minimal number of
projections is simply deduced.  

With the variety of lines intersecting $X$ another problem is also handled. Consider a
set of $N$ points in $\PP^3$ and let $X$ be the curve generated by
these points. Each of these points is projected on a different plane
by a different projection operator $\pi_i$. Those projection operators
are known. We want to recover $X$ by linear computations. Let $Z$ be
the curve, of minimal degree, generated by the centers of projections. We give a formula
of the number of constraints obtained on $X$ as a function of the
degree of $Z$.

Finally we show how those questions were motivated by some problems
related to static and dynamic computational vision. Therefore we
conclude by showing how our results can be applied in that context.

Since our computations will occur in $\PP^3$, we fix $[X,Y,Z,T]^T$, as
homogeneous coordinates, and $T=0$ as the plane at infinity. 

%=============================================
%=============================================
\section{Projection operators}
\label{sec::projMaps}

Let $\pi$ be a projection operator from $\PP^3$ to an embedded
projective plane $i(\PP^2)$ through a point ${\bf O}$. This projection
can be presented by a $3 \times 4$ matrix ${\bf M}$. There exists a
set of simple, but very useful, properties. The kernel of ${\bf M}$ is
exactly the center of projection. The transpose of ${\bf M}$ maps a
line in $i(\PP^2)$ to the plane it defines with the center of
projections, given as point of the dual space $\PP^{3*}$. This can
easily be deduced by a duality argument and a 
simple computation. 

There exists a matrix $\what{\bf M}$, being a polynomial function of
${\bf M}$, which maps a point in $i(\PP^2)$ to the Pl\"ucker
coordinates of the line it generates with the center of
projections. If the matrix ${\bf M}$ is decomposed as follows:
$$
{\bf M}= \left[\begin{array}{c}
	\Ga^T \\
	\La^T \\
	\Th^T
\end{array}\right],
$$
then for ${\bf p}=[x,y,z]^T$, the line ${\bf L}_{\bf p} = \what{{\bf
M}}{\bf p}$ is given by the extensor: ${\bf L}_{{\bf p}}=x\La \meet
\Th + y\Th \meet \Ga + z\Ga \meet \La$, where $\meet$ denotes the meet
operator in the Grassman-Cayley algebra (see
\cite{Barnabei-Brini-Rota-84}). By duality, the matrix $\wtilde{\bf
M}=\what{\bf M}^T$ maps lines in $\PP^3$ to lines in $i(\PP^2)$. 

Consider now two projection operators $\pi_1$ and $\pi_2$. Let ${\bf
O}_1$ and ${\bf O}_2$ be the center of projections and $i_1(\PP^2)$
and $i_2(\PP^2)$ the plane of projections. Let ${\bf e}_j$ be the
point of intersection of $i_j(\PP^2)$ with the line $\overline{{\bf
O}_1 {\bf O}_2}$. Let $\si({\bf e}_2)$ be the pencil of lines in
$i_2(\PP^2)$ through ${\bf e}_2$. It is easy to define a map from
$i_1(\PP^2) \setminus \{{\bf e}_1\}$ to $\si({\bf e}_2)$ as follows. Each
point ${\bf p}$ is sent to the line given by $\pi_2(\pi_1^{-1}({\bf
p}))$. This map is linear and its matrix is ${\bf F} = \wtilde{\bf
M}_2 \what{\bf M}_1$. Following the standard terminology used in
computational vision, we will call the matrix ${\bf F}$ the {\it
fundamental matrix} of the pair of projections $\pi_1$ and
$\pi_2$ and the points ${\bf e}_1$ and ${\bf e}_2$ will be
respectively called the {\it first and the second epipole}. The line
in the first (second) projection plane passing the first (second)
epipole atre called the {\it epipolar lines}. Clearly ${\bf Fe}_1 =
{\bf e}_2^T {\bf F}= {\bf 0}$.

\begin{prop}
The knowledge of the fundamental matrix ${\bf F}$ of a couple of
projections operators allows the recovery of the matrices ${\bf M}_1$
and ${\bf M}_2$ of the projections modulo the action $\tht_2$. More
precisely the couple $({\bf M}_1, {\bf M}_2)$ is equivalent to $([{\bf
I},{\bf O}],[{\bf H},{\bf e}_2])$, where ${\bf H} = -\frac{[{\bf
e}_2]}{\parallel {\bf e}_2 \parallel^2} {\bf
F}$. The matrix  $[{\bf e}_2]$ is defined as being $\tau({\bf e}_2)$
where $\tau$ maps any vector ${\bf x}$ of $\C^3$ to the matrix that represents
the cross-product by ${\bf x}$. We have:
$$
\tau({\bf x}) = [{\bf x}] = \left[\begin{array}{ccc}
	0   & - {\bf x}_3  &  {\bf x}_2 \\
	{\bf x}_3   & 0    &  -{\bf x}_1 \\
	-{\bf x}_2  & {\bf x}_1    & 0
\end{array}\right],
$$
and ${\bf x} = [{\bf x}_1,{\bf x}_2,{\bf x}_3]$.
\end{prop}
\begin{proof}
We are looking for a matrix ${\bf A} \in Pr_3$ such that: 
$$
\begin{array}{ccc}
{\bf M}_1 & = & [{\bf I},{\bf 0}] {\bf A}^{-1} \\
{\bf M}_2 & = & [{\bf H},{\bf e}_2] {\bf A}^{-1} 
\end{array}
$$
Let us write ${\bf B} = {\bf A}^{-1}$ as follows:
$$
{\bf B} = \left[\begin{array}{ccc}
		 {\bf \Om} & \vdots & {\bf u} \\
                 \hdots & \hdots & \hdots \\
		 {\bf v}^T & \vdots & 0
	\end{array}\right]
$$
Let us write ${\bf M}_i = [\overline{\bf M}_i, {\bf m}_i]$. Then it
follows immediately from the definition of ${\bf F}$ that ${\bf F} =
[{\bf e}_2] \overline{\bf M}_2 \overline{\bf M}_1^{-1}$. Using the
following algebraic identity for any two vectors ${\bf x}$ and ${\bf
y}$ in $\C^3$: $({\bf x}^T {\bf y}) {\bf I} = {\bf x}{\bf y}^T - [{\bf
x}][{\bf y}]$, it is easy to prove that it is 
sufficient to take: ${\bf \Om} =\overline{\bf M}_1$, ${\bf u} = {\bf
m}_1$, ${\bf v} = \frac{1}{\parallel {\bf e}_2 \parallel^2}
\overline{\bf M}_2 {\bf e}_2$ and $\lambda = \frac{1}{\parallel {\bf
e}_2 \parallel^2} ({\bf e}_2^T {\bf m}_2 - {\bf e}_2^T {\bf Hm}_1)$.  
\end{proof}

This shows that in order to characterize a couple of projections
modulo the action $\tht_2$, we only need to compute the fundamental
matrix. In the sequel we show how to recover ${\bf F}$ from two
projections of an algebraic curve. 

%==========================================================
%==========================================================
\section{Two projections with unknown projection operators}
\label{sec::twoProj}

In this section we deal with the first problem. A smooth and
irreducible curve $X$ embedded in $\PP^3$ is projected onto two
generic planes through two generic points. The projection operators
$\pi_1$ and $\pi_2$ are unknown. First we want to recover the fundamental
matrix  of the couple $(\pi_1,\pi_2)$ from the projected curve $Y_1 =
\pi_1(X)$ and $Y_2 = \pi_2(X)$.

%=============================
\subsection{Single projection}
\label{subsec:SingleView}

We shall mention of set of well known facts about generic projection.
Let $X$ be a smooth irreducible algebraic embedded in $\PP^3$ and $Y$
its projection on a generic plane
through a center of projections, ${\bf O}$. 
\begin{enumerate}
\item The curve $Y$ will always contain singularities. Furthermore for
a generic position of the center of projection, the only singularities
of $Y$ will be nodes. 
\item The {\it class} of a planar curve is defined to be the degree of
its dual curve. Let $m$ be the class of $Y$. Then $m$ is constant for
a generic position of the center of projections. 
\item If $d$ and $g$ are the degree and the genus of $X$, they are
respectively, for a generic position of ${\bf O}$, the degree and the
genus of $Y$, and the Pl\"ucker formula yields: 
$$
\begin{gathered}
m=d(d-1)-2(\sharp\nodes), \\ g=\frac{(d-1)(d-2)}{2} - (\sharp\nodes),
\end{gathered}
$$
where $\sharp\nodes$ denotes the number of nodes of $Y$. Hence the
genus, the degree and the class are related by 
$$
m=2d+2g-2.
$$
\end{enumerate}

%==============================================
\subsection{Fundamental matrix construction}

We are ready now to investigate the recovery of the fundamental matrix of a
couple of projections $(\pi_1,\pi_2)$ when the only knowledge is made
of the projections of a smooth irreducible curve.

As before, let $X$ be a smooth irreducible
curve embedded in $\PP^3$, which cannot be embedded in a plane. The
degree of $X$ is $d \geq 3$. Let ${\bf M}_i$, $i=1,2$, be the projection
matrices. Let ${\bf F}$, ${\bf e}_1$ and ${\bf e}_2$ be defined as
before. We will need to consider the 
two following mappings: ${\bf p} \stackrel{\ga}{\mapsto} {\bf e}_1
\join {\bf p}$ and ${\bf p} \stackrel{\xi}{\mapsto} {\bf Fp}$, where
$\join$ is the join operator \cite{Barnabei-Brini-Rota-84}, which is
equivalent to the cross-product in that case. Both maps are defined on
the first projection plane.

Let $Y_1$ and $Y_2$ be the projected curves. Assume that they are defined by the polynomials $f_1$
and $f_2$. Let $Y_1^\star$ and $Y_2^\star$ denote the dual 
curves, whose polynomials are respectively $\phi_1$ and $\phi_2$. 

\begin{theorem} \label{theoKrupeq}  
For a generic position of the centers of projections with respect to the curve
$X$, there exists a non-zero scalar $\lambda$, such that for all
points ${\bf p}$ in the projection plane, the following equality holds:  
\begin{equation}\label{ExtKrupEq} 
\phi_2(\xi({\bf p})) = \lambda \phi_1(\gamma({\bf p}))
\end{equation}
For reasons that will be clear later, we shall call this equation, the {\it
generalized Kruppa equation}.
\end{theorem} 
\begin{proof} 
Let $\ep_i$ be the set of epipolar lines tangent to curve in image
$i$. We start with the following lemma: 
\begin{lemma}
The two sets $\ep_1$ and $\ep_2$ are projectively equivalent. Moreover
for each corresponding pair of epipolar lines $({\bf l},{\bf l}') \in
\ep_1 \times \ep_2$, the multiplicities of ${\bf l}$ and ${\bf l}'$ as
points of $Y_1^*$ and $Y_2^*$ are the same.  
\end{lemma}
\begin{proof}
Consider the following three pencils:
\begin{itemize}
\item $\si({\bf L}) \approx \PP^1$, the pencil of {\it epipolar
planes}, that is planes containing
the baseline joining the two centers of projection, 
\item $\si({\bf e}_1) \approx \PP^1$, the pencil of epipolar lines in
the first projection plane, 
\item $\si({\bf e}_2) \approx \PP^1$, the pencil of epipolar lines in
the second projection plane. 
\end{itemize}
Thus we have $\ep_i \subset \si({\bf e}_i)$. Moreover if $E$ is the
set of planes in $\si({\bf L})$ tangent to the curve in space, then
there exists a one-to-one mapping between $E$ and each $\ep_i$ which
leaves the multiplicities unchanged. This completes the proof. 
\end{proof}

This lemma implies that both sides of equation \ref{ExtKrupEq} define
the same algebraic set, that is the the union of epipolar lines
tangent to $Y_1$. Since $\phi_1$ and $\phi_2$, in the generic case,
have the same degree (as stated in \ref{subsec:SingleView}), each side
can be factorized as follows: 
$$
\begin{gathered}
\phi_1(\gamma (x,y,z)) = \prod_i (\alpha_{1i} x + \alpha_{2i} y +
\alpha_{3i} z)^{a_i} \\ \phi_2(\xi (x,y,z)) = \prod_i \lambda_i
(\alpha_{1i} x + \alpha_{2i} y + \alpha_{3i} z)^{b_i},
\end{gathered} 
$$
where $\sum_i a_i=\sum_j b_j=m$. By the previous lemme we also have:
$a_i=b_i$ for all $i$. 
\end{proof}

By eliminating the scalar $\lambda$ from the generalized Kruppa
equation (\ref{ExtKrupEq}) we obtain a set of 
bi-homogeneous equations in ${\bf F}$ and ${\bf e}_1$. Hence they
define a variety in $\PP^2 \times \PP^8$. We turn our attention to the
dimensional analysis of this variety. Our concern is to exhibit the
conditions for which this variety is discrete. From a practical point
of view, this is a step toward the recovery of the original curve.

%==============================
\subsection{Dimension analysis}
\label{sec::dimensionAnalysis}

Let $\{E_i({\bf F}, {\bf e}_1)\}_i$ be the set of bi-homogeneous
equations on ${\bf F}$ and ${\bf e}_1$, extracted from the generalized
Kruppa equation (\ref{ExtKrupEq}). Our first
concern is to determine whether all solutions of equation
(\ref{ExtKrupEq}) are admissible, that is whether they satisfy the
usual constraint ${\bf F}{\bf e}_1=0$. Indeed we prove the following
statement: 

\begin{prop}
As long as there are at least 2 distinct lines through ${\bf e}_1$
tangent to $Y_1$, equation (\ref{ExtKrupEq}) implies that $\rank{\bf
F}=2$ and ${\bf F}{\bf e}_1={\bf 0}$. 
\end{prop}
\begin{proof}
The variety defined by $\phi_1(\gamma({\bf p}))$ is then a union of at
least 2 distinct lines through ${\bf e}_1$.  If equation
(\ref{ExtKrupEq}) holds, $\phi_2(\xi({\bf p}))$ must define the same
variety. 

There are 2 cases to exclude: If $\rank{\bf F}=3$, then the curve
defined by $\phi_2(\xi({\bf p}))$ is projectively equivalent to the
curve defined by $\phi_2$, which is $Y_2^\star$. In particular, it is
irreducible. 

If $\rank{\bf F}<2$ or $\rank{\bf F}=2$ and ${\bf F}{\bf e}_1\neq{\bf
0}$, then there is some ${\bf a}$, not a multiple of ${\bf e}_1$, such
that ${\bf Fa}={\bf 0}$.  Then the variety defined by $\phi_2(\xi({\bf
p}))$ is a union of lines through ${\bf a}$. In neither case can this
variety contain two distinct lines through 
${\bf e}_1$, so we must have $\rank{\bf F}=2$ and ${\bf F}{\bf
e}_1={\bf 0}$. 
\end{proof}

As a result, in a generic situation every solution of $\{E_i({\bf F},
{\bf e}_1)\}_i$ is admissible.  Let $V$ be the subvariety of $\PP^2
\times \PP^8 \times \PP^2$ defined by the equations $\{E_i({\bf
F},{\bf e}_1)\}_i$ together with ${\bf F}{\bf e}_1={\bf 0}$ and ${\bf
e}_2{}^T {\bf F}={\bf 0}^T$, where ${\bf e}_2$ is the second epipole.
We next compute a lower bound on the dimension of $V$, after which we
will be ready for the calculation itself. 

\begin{prop} \label{prop:lowerbound}
If $V$ is non-empty, the dimension of $V$ is at least $7-m$. 
\end{prop}
\begin{proof}
Choose any line ${\bf l}$ in $\PP^2$ and restrict ${\bf e}_1$ to the
affine piece $\PP^2\setminus{\bf l}$. Let $(x,y)$ be homogeneous
coordinates on ${\bf l}$.  If ${\bf F}{\bf e}_1={\bf 0}$, the two
sides of equation (\ref{ExtKrupEq}) are both unchanged by replacing
${\bf p}$ by ${\bf p}+\alpha{\bf e}_1$.  So equation (\ref{ExtKrupEq})
will hold for all ${\bf p}$ if it holds for all ${\bf p}\in{\bf
l}$. Therefore equation (\ref{ExtKrupEq}) is equivalent to the
equality of 2 homogeneous polynomials of degree $m$ in $x$ and $y$,
which in turn is equivalent to the equality of $(m+1)$ coefficients.
After eliminating $\lambda$, we have $m$ algebraic conditions on
$({\bf e}_1,{\bf F},{\bf e}_2)$ in addition to ${\bf Fe}_1={\bf 0}$,
${\bf e}_2{}^T{\bf F}={\bf 0}^T$. 

The space of all epipolar geometries, that is, solutions to ${\bf
F}{\bf e}_1={\bf 0}$, ${\bf e}_2{}^T{\bf F}={\bf 0}^T$, is irreducible
of dimension 7.  Therefore, $V$ is at least $(7-m)$-dimensional. 
\end{proof}

For the calculation of the dimension of $V$ we introduce some
additional notations.  
Given a triplet $({\bf e}_1,{\bf F},{\bf e}_2) \in
\PP^2\times\PP^8\times\PP^2$, let $\{ {\bf 
q}_{1\alpha}({\bf e}_1) \}$ (respectively $\{ {\bf q}_{2\alpha}({\bf
e}_2) \}$) be the tangency points of the epipolar lines through ${\bf
e}_1$ (respectively ${\bf e}_2$) to the first (respectively second)
projected curve. Let ${\bf Q}_\alpha({\bf e}_1,{\bf e}_2)$ be the 3D
points projected onto $\{ {\bf q}_{1\alpha}({\bf e}_1) \}$ and $\{{\bf
q}_{2\alpha}({\bf e}_2) \}$. Let ${\bf L}$ be the baseline joining the
two centers of projections. We next provide  a sufficient condition for $V$ to
be discrete. 

\begin{prop}
For a generic position of the centers of projection, the variety $V$ will be
discrete if, for any point $({\bf e}_1,{\bf F},{\bf e}_2) \in V$, the
union of ${\bf L}$ and the points ${\bf Q}_\alpha({\bf e}_1,{\bf
e}_2)$ is not contained in any quadric surface. 
\end{prop}
\begin{proof}
For generic projections, there will be $m$ distinct points
$\{{\bf q}_{1\alpha}({\bf e}_1)\}$ and $\{{\bf q}_{2\alpha}({\bf
e}_2)\}$, and we can regard ${\bf q}_{1\alpha}$, ${\bf q}_{2\alpha}$
locally as smooth functions of ${\bf e}_1$, ${\bf e}_2$. 

We let $W$ be the affine variety in $\C^3\times\C^9\times\C^3$ defined
by the same equations as $V$.  Let $\Theta=({\bf e}_1, {\bf F}, {\bf
e}_2)$ be a point of $W$ corresponding to a non-isolated point of
$V$. Then there is a tangent vector $\vartheta=({\bf v}, \Phi, {\bf
v}')$ to $W$ at $\Theta$ with $\Phi$ not a multiple of ${\bf F}$. 

If $\chi$ is a function on $W$, $\nabla_{\Theta,\vartheta}(\chi)$ will
denote the derivative of $\chi$ in the direction defined by
$\vartheta$ at $\Theta$.  For 
$$
\chi_\alpha({\bf e}_1,{\bf F},{\bf e}_2)= {\bf q}_{2\alpha}({\bf
e}_2)^T{\bf F}{\bf q}_{1\alpha}({\bf e}_1), 
$$
the generalized Kruppa equation implies that $\chi_\alpha$ vanishes
identically on $W$, so its derivative must also vanish.  This yields 
\begin{equation}\label{EqNablaQFQ}
\begin{split}
\nabla_{\Theta,\vartheta}(\chi_\alpha)&=
(\nabla_{\Theta,\vartheta}({\bf q}_{2\alpha}))^T{\bf F}{\bf
q}_{1\alpha}\\ &\quad+{\bf q}_{2\alpha}^T\Phi{\bf q}_{1\alpha}+ {\bf
q}_{2\alpha}^T{\bf F}(\nabla_{\Theta,\vartheta}({\bf q}_{1\alpha}))\\
&=0.
\end{split}
\end{equation}
We shall prove that $\nabla_{\Theta,\vartheta}({\bf q}_{1\alpha})$ is
in the linear span of ${\bf q}_{1\alpha}$ and ${\bf e}_1$.  (This
means that when the epipole moves slightly, ${\bf q}_{1\alpha}$ moves
along the epipolar line.)  Consider $\kappa(t)=f({\bf
q}_{1\alpha}({\bf e}_1 + t{\bf v}))$, where $f$ is the polynomial
defining the image curve $Y_1$.  Since ${\bf q}_{1\alpha}({\bf e}_1 +
t{\bf v}) \in Y_1$, $\kappa\equiv 0$, so the derivative
$\kappa'(0)=0$.  On the other hand,
$\kappa'(0)=\nabla_{\Theta,\vartheta}(f({\bf q}_{1\alpha}))
=\grad_{{\bf q}_{1\alpha}}(f)^T\nabla_{\Theta,\vartheta}({\bf
q}_{1\alpha})$. 

Thus we have $\grad_{{\bf
q}_{1\alpha}}(f)^T\nabla_{\Theta,\vartheta}({\bf q}_{1\alpha})=0$. But
also $\grad_{{\bf q}_{1\alpha}}(f)^T{\bf q}_{1\alpha}=0$ and
$\grad_{{\bf q}_{1\alpha}}(f)^T{\bf e}_1=0$.  Since $\grad_{{\bf
q}_{1\alpha}}(f)\neq{\bf O}$ (${\bf q}_{1\alpha}$ is not a singular
point of the curve), this shows that $\nabla_{\Theta,\vartheta}({\bf
q}_{1\alpha})$, ${\bf q}_{1\alpha}$, and ${\bf e}_1$ are linearly
dependent.  ${\bf q}_{1\alpha}$ and ${\bf e}_1$ are linearly
independent, so $\nabla_{\Theta,\vartheta}({\bf q}_{1\alpha})$ must be
in their linear span. 

We have that  ${\bf q}_{2\alpha}^T{\bf F}{\bf e}_1= {\bf
q}_{2\alpha}^T{\bf Fq}_{1\alpha}=0$, so ${\bf q}_{2\alpha}^T{\bf F}
\nabla_{\Theta,\vartheta}({\bf q}_{1\alpha})=0$: the third term of
equation (\ref{EqNablaQFQ}) vanishes.  In a similar way, the first
term of equation (\ref{EqNablaQFQ}) vanishes, leaving 
\begin{equation}
\label{q_al}
{\bf q}_{2\alpha}^T\Phi{\bf q}_{1\alpha}=0. 
\end{equation}

The derivative of $\chi({\bf e}_1,{\bf F},{\bf e}_2)={\bf F}{\bf e}_1$
must also vanish, which  yields 
\begin{equation}
\label{e_12}
{\bf e}_2{}^T\Phi{\bf e}_1=0.
\end{equation}
   
From equality (\ref{q_al}), we deduce that for every ${\bf Q}_\alpha$,
we have 
$$
{\bf Q}^T_\alpha {\bf M}_2^T \Phi {\bf M}_1 {\bf Q}_\alpha=0.
$$
From equality (\ref{e_12}), we deduce that every point ${\bf P}$ lying
on the baseline must satisfy 
$$
{\bf P}^T {\bf M}_2^T \Phi {\bf M}_1 {\bf P}=0.
$$
The fact that $\Phi$ is not a multiple of $F$ implies that ${\bf
M}_2^T\Phi{\bf M}_1 \neq 0$, so together these two last equations mean
that the union ${\bf L} \cup \{{\bf Q}_\alpha\}$ lies on a quadric 
surface.  Thus if there is no such quadric surface, every point in $V$
must be isolated. 
\end{proof}

Observe that this result is consistent with the previous proposition,
since there always exist a quadric surface containing a given line and
six given points. However in general there is no quadric containing a 
given line and seven given points. Therefore we can conclude with the
following theorem. 

\begin{theorem}
For a generic position of the centers of projection, the generalized Kruppa
equation defines the epipolar geometry up to a finite-fold ambiguity
if and only if $m \geq 7$. 
\end{theorem}

Since different curves in generic position give rise to independent
equations, this result means that the sum of the classes of the projected
curves must be at least $7$ for $V$ to be a finite set.  

%===============================
\subsection{Recovering the curve}

Let the projection matrices be ${\bf M}_1$ and ${\bf
M}_2$. Hence the two cones defined by the projected curves and the
centers of projections are given by: $\Delta_1({\bf P})=f_1({\bf M}_1 {\bf P})$ and
$\Delta_2({\bf P})=f_2({\bf M}_2 {\bf P})$. The reconstruction is
defined as the curve whose equations are $\Delta_1=0$ and
$\Delta_2=0$. This curve has two irreducible components as the
following theorem states. 

\begin{theorem}
\label{Theo::TwoViews3Drecons}
For a generic position of the centers of projection, namely when no epipolar
plane is tangent twice to the curve $X$, the curve defined by
$\{\Delta_1=0,\Delta_2=0\}$ has two irreducible components. One has
degree $d$ and is the actual solution of the reconstruction. The other
one has degree $d(d-1)$. 
\end{theorem}
\begin{proof}
For a line ${\bf L}\subset\PP^3$, we write $\sigma({\bf L})$ for the
pencil of planes containing ${\bf L}$.  For a point ${\bf p}\in\PP^2$,
we write $\sigma({\bf p})$ for the pencil of lines through ${\bf
p}$. There is a natural isomorphism between $\sigma({\bf e}_i)$, the
epipolar lines in image $i$, and $\sigma({\bf L})$, the planes
containing both centers of projections.  Consider the following covers of
$\PP^1$: 
\begin{enumerate}
\item $X \stackrel{\eta}{\longrightarrow} \sigma({\bf L}) \cong
\PP^1$, taking a point $x \in X$ to the epipolar plane that it defines
with the centers of projection. 
\item $Y_1 \stackrel{\eta_1}{\longrightarrow} \sigma({\bf e}_1) \cong
\sigma({\bf L}) \cong \PP^1$, taking a point $y \in Y_1$ to its
epipolar line in the first projection plane. 
\item $Y_2 \stackrel{\eta_2}{\longrightarrow} \sigma({\bf e}_2) \cong
\sigma({\bf L}) \cong \PP^1$, taking a point $y \in Y_2$ to its
epipolar line in the second projection plane. 
\end{enumerate}
If $\rho_i$ is the projection $X \to Y_i$, then
$\eta=\eta_i\rho_i$. Let ${\cal B}$ the union set of branch points of
$\eta_1$ and $\eta_2$.  It is clear that the branch points of $\eta$
are included in ${\cal B}$.  Let $S=\PP^1 \setminus {\cal B}$, pick $t
\in S$, and write $X_S=\eta^{-1}(S)$, $X_t=\eta^{-1}(t)$.  Let
$\mu_{X_S}$ be the monodromy: $\pi_1(S,t) \longrightarrow \Perm(X_t)$,
where $\Perm(Z)$ is the group of permutation of a finite set $Z$, see
\cite{Miranda-91}. It is well known that the path-connected components
of $X$ are in one-to-one correspondence with the orbits of the action of
$\im(\mu_{X_S})$ on $X_t$.  Since $X$ is assumed to be irreducible, it
has only one component and $\im(\mu_{X_S})$ acts transitively on
$X_t$.  Then if $\im(\mu_{X_S})$ is generated by transpositions, this
will imply that $\im(\mu_{X_S})=\Perm(X_t)$.  In order to show that
$\im(\mu_{X_S})$ is actually generated by transpositions, consider a
loop in $\PP^1$ based at $t$, say $l_t$.  If $l_t$ does not go round
any branch point, then $l_t$ is homotopic to the constant path in $S$
and then $\mu_{X_S}([l_t])=1$.  Now in ${\cal B}$, there are three
types of branch points: 
\begin{enumerate}
\item branch points that come from nodes of $Y_1$: these are not
branch points of $\eta$, 
\item branch points that come from nodes of $Y_2$: these are not
branch points of $\eta$, 
\item branch points that come from epipolar lines tangent either to
$Y_1$ or to $Y_2$: these are genuine branch points of $\eta$. 
\end{enumerate}   

If the loop $l_t$ goes round a point of the first two types, then it
is still true that $\mu_{X_S}([l_t])=1$.  Now suppose that $l_t$ goes
round a genuine branch point of $\eta$, say $b$ (and goes round no
other points in $\cal B$).  By genericity, $b$ is a simple two-fold
branch point, hence $\mu_{X_S}([l_t])$ is a transposition.  This shows
that $\im(\mu_{X_S})$ is actually generated by transpositions and so
$\im(\mu_{X_S})=\Perm(X_t)$. 

Now consider $\tilde{X}$, the curve defined by
$\{\Delta_1=0,\Delta_2=0\}$.  By Bezout's Theorem $\tilde{X}$ has
degree $d^2$.  Let $\tilde{x} \in \tilde{X}$.  It is projected onto a
point $y_i$ in $Y_i$, such that $\eta_1(y_1)=\eta_2(y_2)$.  Hence
$\tilde{X} \cong Y_1 \times_{\PP^1} Y_2$; restricting to the inverse
image of the set $S$, we have $\tilde{X}_S \cong X_S \times_S X_S$. We
can therefore identify $\tilde{X}_t$ with $X_t \times X_t$.  The
monodromy $\mu_{\tilde{X}_S}$ can then be given by
$\mu_{\tilde{X}_S}(x,y)=(\mu_{X_S}(x),\mu_{X_S}(y))$.  Since
$\im(\mu_{X_S})=\Perm(X_t)$, the action of $\im(\mu_{\tilde{X}_S})$ on
$X_t \times X_t$ has two orbits, namely $\{(x,x)\} \cong X_t$ and
$\{(x,y) | x \neq y\}$.  Hence $\tilde{X}$ has two irreducible
components.  One has degree $d$ and is $X$, the other has degree
$d^2-d=d(d-1)$. 
\end{proof} 
  
This result provides a way to find the right solution for the
recovery in a generic configuration, except in the case of
conics, where the two components of the reconstruction are both
admissible. 

%=========================================================================
%=========================================================================
\section{The $N >> 1$ projections problem with known projection operators}
\label{sec::Nproj}

Now we turn our attention to the second problem. $N >> 1$ projections
maps $\{\pi_i\}_{i=1,...,N}$ given by $N$ matrices $\{{\bf
M}_i\}_{i=1,...,N}$ are known. Therefore $N$ projections of an
irreducible smooth algebraic curve are also provided. The problem is
to recover the original curve by linear computations as much as possible.

%================================================
\subsection{Curve presentation in the dual space}
\label{sec::Dual}

Let $X^\star$ be the dual variety of $X$. Since $X$ is supposed not to
be a line, the dual variety $X^\star$ must be a hypersurface of the
dual space \cite{Harris-92}. Hence let ${\Upsilon}$ be a minimal
degree polynomial that represents $X^\star$. Our first concern is to
determine the degree of $\Upsilon$. 

\begin{prop}
The degree of $\Upsilon$ is $m$, that is, the common degree of the
dual projected curves. 
\end{prop}
\begin{proof}
Since $X^\star$ is a hypersurface of $\PP^{3\star}$, its degree is the
number of points where a generic line in $\PP^{3\star}$ meets
$X^\star$. By duality it is the number of planes in a generic pencil
that are tangent to $X$. Hence it is the degree of the dual projected
curve. Another way to express the same fact is the observation that
the dual projected curve is the intersection of $X^\star$ with a generic
plane in $\PP^{3\star}$. Note that this provides a new proof that the
degree of the dual projected curve is constant for a generic position of
the center of projection. 
\end{proof}

For the recovery of $X^*$ from multiple projection, we will need to
consider the mapping from a line ${\bf l}$ of the projection plane to the
plane that it defines with the center of projection. Let $\mu: {\bf l}
\mapsto {\bf M}^T {\bf l}$ denote this mapping. There exists a link
involving $\Upsilon$, $\mu$ and $\phi$, the polynomial of the dual
projected curve: $\Upsilon(\mu({\bf l}))=0$ whenever $\phi({\bf l})=0$.
Since these two polynomials have the same degree (because $\mu$ is
linear) and $\phi$ is irreducible, there exists a scalar $\lambda$
such that 
$$
\Upsilon(\mu({\bf l}))=\lambda \phi({\bf l})
$$
for all lines ${\bf l} \in \PP^{2\star}$. Eliminating $\lambda$, we
get $\binom{m+2}{m}-1$ linear equations on $\Upsilon$.  Since the
number of coefficients in $\Upsilon$ is $\binom{m+3}{m}$, we can state
the following result: 

\begin{prop} \label{dualNbViews}
The recovery in the dual space can be done linearly using at
least $k \geq \frac{m^2+6m+11}{3(m+3)}$ projections. 
\end{prop}

%======================================================================
\subsection{Curve presentation in the Grasmmannian of lines of $\PP^3$}
\label{sec::Grass}

Let $\G(1,3)$ be the Grassmanian of lines of $\PP^3$. Consider the set
of lines in $\PP^3$ intersecting the curve $X$ of degree $d$. This defines
a subvariety of $\G(1,3)$ which is the intersection of $\G(1,3)$ with
a hypersurface of degree $d$ in $\PP^5$, given by a homogeneous
polynomial $\Gamma$, defined modulo the $d$th graded piece
$I(\G(1,3))_d$ of the ideal of $\G(1,3)$ and modulo scalars. However
picking one representative of this equivalence class is sufficient to
recover entirely without any ambiguity the curve $X$. In our
context, we shall call any representative of this class the {\it Chow
polynomial} of the curve. We need to compute the class of $\Ga$ in the
homogeneous coordinate ring of $\G(1,3)$, or more precisely in its
$d$th graded piece, $S(\G(1,3))_d$, which dimension is $Nd =
\binom{d+5}{d} - \binom{d-2+5}{d-2}$.  

Let $f$ be the polynomial defining the projected curve, $Y$. Consider the
mapping  that associates to a point in the projection plane the line
it generates with the center of projection: $\nu: {\bf
p} \mapsto \widehat{\bf M} {\bf p}$. The polynomial $\Gamma(\nu({\bf
p}))$ vanishes whenever $f({\bf p})$ does. Since they have same degree
and $f$ is irreducible, there exists a scalar $\lambda$ such that for
every point ${\bf p} \in \PP^2$, we have:   
$$
\Gamma(\nu({\bf p})) = \lambda f({\bf p}).
$$
This yields $\binom{d+2}{d}-1$ linear equations on $\Gamma$.  

Hence a similar statement to that in Proposition \ref{dualNbViews} can
be made: 
\begin{prop}
The recovery in $\G(1,3)$ can be done linearly using at least $k
\geq \frac{1}{6}\frac{d^3+5d^2+8d+4}{d}$ projections. 
\end{prop}

%==============================================================================
\subsection{Family of projections operators and finite closed subset of points} 

Consider now a finite collection of points ${\bf P}_i$ in $\PP^3$. Each point is
projected by a different projection map. The $N$ maps are known and so
the projected points. 

Let $X$ be the smooth irreducible curve generated by the points ${\bf
P}_i$ and $Y$ be the smooth irreducible curve, of minimal degree, generated by the center
of projections. 

Each projected point ${\bf p}_i$ yiels one linear equation on the
variety of intersecting lines of $X$, namely $\Ga(\pi_i^{-1}({\bf
p}_i)) = 0$, where $\Ga$ is the Chow polynomial of $X$ as before. 

Let ${\bf d}$ and ${\bf d}'$ be respectively the degree of $X$ and
$Y$. We compute the number of constraints obtained on $\Ga$ from the
projected points as a function of ${\bf d}$ and ${\bf d}'$. In other
words, we want to compute the maximal number of constraints that one
can extract on a smooth irreducible curve $X$ embedded in $\PP^3$ from
a finite number of lines in the join of $X$ with a known curve $Y$.

\begin{prop}
\label{Theo::numbOfConst}
The maximal number of constraints is $N_d - (h^0({\cal
O}_{\PP^5}(d-d')) - h^0({\cal O}_{\PP^5}(d-d'-2))+1)$
\footnote{As usual, $h^0({\cal O}_{\PP^5}(k))$ denotes the dimension
of the cohomology group $H^0(\PP^5,{\cal O}_{\PP^5}(k))$ (see
\cite{Hartshorne-77}).}, 
where $N_d = dim(S(\G(1,3))_d)$ is the dimension of the $d$-th graded piece of
homogeneous coordinate ring of $\G(1,3)$.  
\end{prop}
\begin{proof}
Each projected point generates a line with the center of projection. 
Let ${\bf L}_1,....,{\bf L}_n$ be these
$n$ lines joining $X$ and $Y$. Let $\Ga_X$ and $\Ga_Y$ be the Chow
polynomial of $X$ and $Y$ respectively. We shall denote by $Z(\Ga_X)$
and $Z(\Ga_Y)$ the sets where they vanish. Let $V = Z(\Ga_X) \cap
Z(\Ga_Y) \cap \G(1,3)$. For $n >> 1$, we have 
$$
\begin{array}{l}
\{ \Ga \in H^0(\PP^5,{\cal O}_{\PP^5}(d)) \,:\, \Ga(L_i)=0, i=1,\dots, n\} = \\ 
\{ \Ga \in H^0(\PP^5,{\cal O}_{\PP^5}(d)) \,:\, \Ga_V \equiv 0 \} =
I_{V,\PP^5}(d). 
\end{array}
$$
So, we want to compute $dim(I_{V,\PP^5}(d))$, or equivalently, $h^0(V,{\cal O}_V(d)) = h^0({\cal
O}_{\PP^5}(d)) - dim(I_{V,\PP^5}(d))$. Since $V$ is a complete
intersection of degree $(d, d', 2)$ in $\PP^5$,  the dimension of
$I_{V,\PP^5}(d)$ should be equal to 
$$
h^0({\cal O}_{\PP^5}(d-2)) + h^0({\cal O}_{\PP^5}(d-d')) - h^0({\cal
O}_{\PP^5}(d-d'-2))+1. 
$$

As a consequence
$$
h^0(V,{\cal O}_V(d)) = N_d - 
\left( 
h^0({\cal O}_{\PP^5}(d-d')) - h^0({\cal O}_{\PP^5}(d-d'-2)) + 1  
\right).
$$
\end{proof}

%================================================================
%================================================================
\section{Applications to static and dynamic computational vision}

The results obtained above were motivated by some applications to
computational vision. We now proceed to show how these results can be
applied to this field. We start by a quick survey on linear
computational vision. More details can be found in
\cite{Faugeras-93,Hartley-Zisserman-00,Faugeras-Luong-01}. Some of the terminology was introduced before in
section \ref{sec::projMaps}.

%======================================================
\subsection{Foundations of linear computational vision}

Projective algebraic geometry provides a natural framework to
geometric computer vision. However one has to keep in mind that the
geometric entities to be considered are in fact embedded in the
physical three-dimensional Euclidean space. Euclidean space is
provided with three structures defined by three groups of
transformations: the orthogonal group $Euc_3$ (which defines the
Euclidean structure and which is included into the affine group),
$Aff_3$ (defining the affine structure and itself included into the
projective group), and $Pr_3$ (defining the projective structure). We
fix $[X,Y,Z,T]^T$, as homogeneous coordinates, and $T=0$ as the plane
at infinity.

%==================================
\subsection{A single camera system}
\label{SingleCamera}

Computational vision starts with images captured by cameras. The
camera components are the following: 
\begin{itemize}
\item a plane ${\cal R}$, called the {\it retinal plane} or {\it image
plane}; 
\item a point ${\bf O}$, called the {\it optical centre} or {\it
camera centre}, which does not lie on the plane ${\cal R}$. 
\end{itemize}

The plane ${\cal R}$ is regarded as a two dimension projective space
embedded into $\PP^3$. Hence it is also denoted by $i(\PP^2)$. The
camera is a projection machine:  $\pi: \PP^3 \setminus \{{\bf O}\}
\rightarrow i(\PP^2), {\bf P} \mapsto \overline{{\bf O} {\bf P}} \cap
i(\PP^2)$. The projection $\pi$ is determined (up to a scalar) by a $3
\times 4$ matrix ${\bf M}$ (the image of ${\bf P}$ being $\la {\bf M
P}$).  

The physical properties of a camera imply that ${\bf M}$ can be
decomposed as follows: 
$$
{\bf M} = \left[\begin{array}{ccc}
	f & s & u_0 \\
	0 & \al f & v_0 \\
	0 & 0 & 1 
	\end{array}\right] [{\bf R}; {\bf t}],
$$
where $(f,\al,s,u_0,v_0)$ are the so-called internal parameters of the
camera, whereas the rotation ${\bf R}$ and the translation ${\bf t}$
are the external parameters. 

It is easy to see that:
\begin{itemize}
\item The camera centre ${\bf O}$ is given by ${\bf MO}={\bf 0}$. 
\item The matrix ${\bf M}^T$ maps a line in $i(\PP^2)$ to the only
plane containing both the line and ${\bf O}$. 
\item There exists a matrix $\what{\bf M} \in {\cal M}_{6 \times
3}(\R)$, which is a polynomial function of ${\bf M}$, that maps a
point ${\bf p} \in i(\PP^2)$ to the line $\overline{{\bf Op}}$
(optical ray), represented by its Pl\"{u}cker coordinates in
$\PP^5$. If the camera matrix is decomposed as follows: 
 $$
{\bf M}= \left[\begin{array}{c}
	\Ga^T \\
	\La^T \\
	\Th^T
\end{array}\right],
$$
then for ${\bf p}=[x,y,z]^T$, the optical ray ${\bf L}_{\bf p} =
\what{{\bf M}}{\bf p}$ is given by the extensor: ${\bf L}_{{\bf
p}}=x\La \meet \Th + y\Th \meet \Ga + z\Ga \meet \La$, where $\meet$
denotes the meet operator in the Grassman-Cayley algebra (see
\cite{Barnabei-Brini-Rota-84}). 
\item The matrix $\wtilde{\bf M}=\what{\bf M}^T$ maps lines in $\PP^3$
to lines in $i(\PP^2)$. 
\end{itemize} 

Moreover we will need in the sequel to consider the projection of the
{\bf absolute conic} onto the image plane. The absolute conic is
simply defined by the following equations: 
$$
\left\{\begin{array}{rcl}
X^2+Y^2+Z^2 & = & 0 \\
T & = & 0
\end{array}\right.
$$

By definition, the absolute conic is left invariant by Euclidean
transformations. Therefore its projection onto the image plane,
defined by the matrix $\om$, is a function of the internal parameters
only. By Cholesky decomposition $\om={\bf L} {\bf U}$, where ${\bf L}$
(respectively ${\bf U}$) is lower (respectively upper) triangular
matrix. Hence it is easy to see that ${\bf U}=\overline{\bf M}^{-1}$,
where $\overline{\bf M}$ is the $3 \times 3$ matrix of the internal
parameters of ${\bf M}$.

%===================================
\subsection{A system of two cameras}

Given two cameras, $({\bf O}_j, i_j(\PP^2))_{j=1,2}$ are their
components where $i_1(\PP^2)$ and $i_2(\PP^2)$ are two generic
projective planes embedded into $\PP^3$, and  ${\bf O}_1$ and ${\bf
O}_2$ are two generic points in $\PP^3$ not lying on the above
planes. As in \ref{SingleCamera}, let $\pi_j: \PP^3 \setminus \{ {\bf
O}_j \} \ra i_j(\PP^2), {\bf P} \mapsto \overline{{\bf O}_j {\bf P}}
\cap i_j(\PP^2)$ be the respective projections. The camera matrices
are ${\bf M}_i,i=1,2$. 

\subsubsection{Homography between two images of the same plane} 

Consider the case where the two cameras are looking at the same plane
in space, denoted by $\De$. Let  
$$
{\bf M}_i = \left[\begin{array}{c}
	\Ga_i^T \\
	\La_i^T \\
	\Th_i^T
\end{array}\right]
$$
be the camera matrices, decomposed as above.  Let ${\bf P}$ be a point
lying on $\De$. We shall denote the projections of ${\bf P}$ by ${\bf
p}_i=[x_i,y_i,z_i]^T \cong {\bf M}_i {\bf P}$, where $\cong$ means
equality modulo multiplication by a non-zero scalar.  

The optical ray generated by ${\bf p}_1$ is given by ${\bf L}_{{\bf
p}_1}=x_1\La_1 \meet \Th_1 + y_1\Th_1 \meet \Ga_1 + z_1\Ga_1 \meet
\La_1$. Hence ${\bf P}={\bf L}_{{\bf p}_1} \meet \De = x_1\La_1 \meet
\Th_1 \meet \De + y_1\Th_1 \meet \Ga_1 \meet \De + z_1\Ga_1 \meet
\La_1 \meet \De$. Hence ${\bf p}_2 \cong {\bf M}_2 {\bf P}$ is given
by the following expression: ${\bf p}_2 \cong {\bf H}_\De {\bf p}_1$
where: 
$$
{\bf H}_\De = \left[ \begin{array}{ccc}  
     \Ga_2^T (\La_1 \meet \Th_1 \meet \De)  & \Ga_2^T (\Th_1 \meet
\Ga_1 \meet \De) & \Ga_2^T (\Ga_1 \meet \La_1 \meet \De) \\ 
     \La_2^T (\La_1 \meet \Th_1 \meet \De)  & \La_2^T (\Th_1 \meet
\Ga_1 \meet \De) & \La_2^T (\Ga_1 \meet \La_1 \meet \De) \\ 
     \Th_2^T (\La_1 \meet \Th_1 \meet \De)  & \Th_2^T (\Th_1 \meet
\Ga_1 \meet \De) & \Th_2^T (\Ga_1 \meet \La_1 \meet \De) \\ 
\end{array}\right].
$$
This yields the expression of the collineation ${\bf H}_\De$ between
two images of the same plane. 

\begin{definition}
The previous collineation is called the {\em homography} between the
two images, through the plane $\De$. 
\end{definition}

\subsubsection{Epipolar geometry}

\begin{definition}
Let $({\bf O}_j, i_j(\PP^2), {\bf M}_j)_{j=1,2}$ be defined as
before. Given a pair $({\bf p}_1,{\bf p}_2) \in i_1(\PP^1) \times
i_2(\PP^2)$, we say that it is a pair of {\em corresponding} or {\em
matching} points if there exists ${\bf P} \in \PP^3$ such that ${\bf
p}_j=\pi_j({\bf P})$ for $j=1,2$.  
\end{definition}

Consider a point ${\bf p} \in i_1(\PP^2)$. Then ${\bf p}$ can be the
image of any point lying on the fiber $\pi_1^{-1}({\bf p})$. The
matching point in the second image must lie on $\pi_2(\pi_1^{-1}({\bf
p}))$, which is, for a generic point ${\bf p}$, a line on the second
image. Since $\pi_1$ and $\pi_2$ are both linear, there exists a
matrix ${\bf F} \in {\cal M}_{3 \times 3}(\R)$, such that: $\xi({\bf
p})=\pi_2(\pi_1^{-1}({\bf p}))={\bf Fp}$ for all but one point in the
first image.  

\begin{definition}
The matrix ${\bf F}$ is called the {\em fundamental matrix}, where as
the line ${\bf l}_{\bf p}={\bf Fp}$ is called the {\em epipolar line}
of ${\bf p}$. 
\end{definition} 

Let ${\bf e}_1=\overline{{\bf O}_1 {\bf O}_2} \cap i_1(\PP^2)$ and
${\bf e}_2=\overline{{\bf O}_1 {\bf O}_2} \cap i_2(\PP^2)$. Those two
points are respectively called the {\it first} and the {\it second
epipole}. It is easy to see that ${\bf Fe}_1={\bf 0}$, since
$\pi_1^{-1}({\bf e}_1)=\overline{{\bf O}_1 {\bf O}_2}$ and
$\pi_2(\overline{{\bf O}_1 {\bf O}_2})={\bf e}_2$. Observe that by
symmetry ${\bf F}^T$ is the fundamental matrix of the reverse couple
of images. Hence ${\bf F}^T {\bf e}_2={\bf 0}$. Since the only point
in the first image that is mapped to zero by ${\bf F}$ is the first
epipole, ${\bf F}$ has rank $2$.  

Now we want to deduce an expression of ${\bf F}$ as a function of the
camera matrices. By the previous analysis, it is clear that ${\bf
F}=\wtilde{\bf M}_2 \what{\bf M}_1$. Moreover we have the following
properties: 

\begin{prop}
For any plane $\De$, not passing through the camera centres, the
following equalities hold: 
\begin{enumerate}
\item $$ {\bf F} \cong [{\bf e}_2]_\times {\bf H}_\De, $$ where $[{\bf
e}_2]_\times$ is the matrix associated with the cross-product as
follows: for any vector ${\bf p}$, ${\bf e}_2 \times {\bf p}=[{\bf
e}_2]_\times {\bf p}$. Hence we have: 
$$
[{\bf e}_2]_\times = \left[\begin{array}{ccc}
	0            & - {{\bf e}_2}_3  &  {{\bf e}_2}_2 \\
	{{\bf e}_2}_3   & 0             &  -{{\bf e}_2}_1 \\
	-{{\bf e}_2}_2  & {{\bf e}_2}_1    & 0
\end{array}\right].$$	
In particular, we have: ${\bf F}=[{\bf e}_2]_\times {\bf H}_\infty$,
where ${\bf H}_\infty$ is the homography between the two images
through the plane at infinity. 
\item 
\begin{eqnarray} \label{FH} 
{\bf H}_\De^T {\bf F} + {\bf F}^T {\bf H}_\De={\bf 0}.
\end{eqnarray}
\end{enumerate}
\end{prop}
\begin{proof}
The first equality is clear according to its geometric meaning. Given
a point ${\bf p}$ in the first image, ${\bf Fp}$ is its epipolar line
in the second image. The optical ray ${\bf L}_{\bf p}$ passing trough
${\bf p}$ meets the plane $\De$ in a point ${\bf Q}$, which projection
in the second image is ${\bf H}_\De {\bf p}$. Hence the epipolar line
must be ${\bf e}_2 \join {\bf H}_\De {\bf p}$. This gives the required
equality. The second equality is simply deduced from the first one by
a short calculation. 
\end{proof}

\begin{prop}
For a generic plane $\De$, the following equality hold: 
$$
{\bf H}_\De {\bf e}_1 \cong {\bf e}_2.
$$
\end{prop}
\begin{proof}
The image of ${\bf e}_1$ by the homography must be the projection on
the second image of the point defined as being the intersection of the
optical ray generated by ${\bf e}_1$ and the plane $\De$. Hence ${\bf
H}_\De {\bf e}_1 = {\bf M}_2 ({\bf L}_{{\bf e}_1} \meet \De)$. But
${\bf L}_{{\bf e}_1} = \overline{{\bf O}_1 {\bf O}_2}$. Thus the
result must be ${\bf M}_2 {\bf O}_1$ (except when the plane is passing
through ${\bf O}_2$) that is the second epipole ${\bf e}_2$. 
\end{proof}

\subsubsection{Canonical stratification of the reconstruction}

Three-dimension reconstruction can be achieved from a system of two
cameras, once the camera matrices are known. However a typical
situation is that the camera matrices are unknown. Then we face a
double problem: recovering the camera matrices and the actual
object. There exists an inherent ambiguity. Consider a pair of camera
matrices $({\bf M}_1, {\bf M}_2)$. If you change the world coordinate
system by a transformation ${\bf V} \in Pr_3$, the camera matrices are
mapped to $({\bf M}_1 {\bf V}^{-1},{\bf M}_2 {\bf V}^{-1})$. Therefore
we define the following equivalence relation: 

\begin{definition}
Given a group of transformations $G$, two pairs of camera matrices,
say $({\bf M}_1,{\bf M}_2)$ and $({\bf N}_1, {\bf N}_2)$, are said to
be equivalent modulo $G$ if there exists ${\bf V} \in G$ such that
${\bf M}_1 = {\bf V} {\bf N}_1$ and ${\bf M}_2 = {\bf V} {\bf N}_2$. 
\end{definition}

Any reconstruction algorithm will always yields a reconstruction
modulo a certain group of transformations. More presicely there exist
three levels of reconstruction according to the information that can
be extracted from the two images and from a priori knowledge of the
world. \\ 

{\bf Projective stratum.} 

When the only available information is the fundamental matrix, then
the reconstruction is done modulo $Pr_3$. Indeed, from ${\bf F}$, the
so-called intrinsic homography ${\bf S}=-\frac{{\bf e}_2}{\parallel
{\bf e}_2 \parallel} {\bf F}$ is computed and the camera matrices are
equivalent to $([{\bf I};{\bf 0}],[{\bf S};{\bf e}_2])$. \\ 

{\bf Affine stratum.} 

When, in addition to the epipolar geometry, the homography between the
two images through the plane at infinity, denoted by ${\bf H}_\infty$,
can be computed, the reconstruction can be done modulo the group of
affine transformations. Then the two camera matrices are equivalent
$([{\bf I};{\bf 0}],[{\bf H}_\infty;{\bf e}_2])$.\\ 

{\bf Euclidean stratum.}

The Euclidean stratum is obtained by the data of the projection of the
absolute conic $\Om$ onto the image planes, which allows the recovery
the internal parameters of the cameras. Once the internal parameters
of the cameras are known, the relative motion between the cameras
expressed by a rotation ${\bf R}$ and a translation ${\bf t}$ can be
extracted from the fundamental matrix. However only the direction of
${\bf t}$, not the norm, can be recovered. Then the cameras matrices
are equivalent, modulo the group of similarity transformations, to
$(\overline{\bf M}_1 [{\bf I};{\bf 0}], \overline{\bf M}_2 [{\bf
R};{\bf t}])$, where $\overline{\bf M}_1$ and $\overline{\bf M}_2$ are
the matrices of internal parameters.  

Note that the projection of the absolute conic on the image can be
computed using some a priori knowledge of the world. Moreover there
exist famous equations linking $\om_1$ and $\om_2$, the two matrices
defining the projection of the absolute conic onto the images, when
the epipolar geometry is given. This is  the so-called {\em Kruppa
equation}, defined in the following proposition. 

\begin{prop}
\label{KrupEq}
The projections of the absolute conic onto two images are related as
follows. There exists a scalar $\la$ such that 
$$
[{\bf e}_1]_\times^T \om_1^* [{\bf e}_1]_\times = \la {\bf F}^T \om_2^* {\bf F},
$$
where $[{\bf e}_1]_\times$ is the matrix representing the
cross-product by ${\bf e}_1$ and $\om_i^*$ is the adjoint matrix of
$\om_i$. 
\end{prop}

Let ${\ep}_i$ be the tangents to $\pi_i(\Om)$ through ${\bf
e}_i$. Kruppa equation simply states that ${\ep}_1$ and ${\ep}_2$
are projectively isomorphic.  

%========================================================================
\subsection{Applications of the previous results to computational vision}

The previous results (sections \ref{sec::twoProj} and \ref{sec::Nproj}) can be applied to
computational vision in different contexts: 
\begin{enumerate}

\item The recovery of the epipolar geometry from two images of the
same smooth irreducible curve. Theorem \ref{theoKrupeq}
generalizes Kruppa equation to algeraic curves. Section
\ref{sec::dimensionAnalysis} provides with a necessary and sufficient
conditions on the degree of the curve for the epipolar geometry to be
defined up to finite fold ambiguity. Note that the case of conic
sections was first introduced in
\cite{Kahl-Heyden-98,Kaminski-Shashua-00}.

\item The 3D reconstruction of a curve from two images is possible in a
generic situation as shown in theorem
\ref{Theo::TwoViews3Drecons}. The case of conic was also treated in
\cite{Kaminski-Shashua-00,Quan-96,Schmid-Zisserman-98}. Note that
\cite{Forsyth} presents an algorithm for curve reconstruction using a 
blow-up of the projected curve. This nice result, however, does not
provide any information about the relative position of the curve in
$\PP^3$ with respect to other elements of the scene. On the other hand,
our approach based on two images allows reconstructing the curve in
the context of the whole scene. Furthermore the problem of curve
reconstruction was also considered in \cite{Berthilsson-all-02} from
the point of view of global optimization and bundle ajustement. Our
approach, on the contrary, is based on looking at algebraic curves for
which the represenation is more compact.

\item The 3D reconstruction of a curve from $N >> 1$ projection is
linear using the dual space or the Grassmannian of lines $\G(1,3)$,
sections \ref{sec::Dual} and \ref{sec::Grass}. The formalism of the
dual space in the case of conic or quadric was also used in
\cite{Cross-Zisserman-98,Ma-Chen-94}. 

\item The trajectory recovery of a moving point viewed by a moving
camera whose matrix is known over time is a linear problem when using
the variety of intersecting lines of the curve generated by the
motion of the point. Moreover this gives rise to the question of
counting the number of constraints that can be obtained. This is done
in theorem \ref{Theo::numbOfConst}. Note that our algorithm for
trajectory recovery or triangulation is a complete generalization of
\cite{Avidan-Shashua-00}. 
\end{enumerate}

%=====================================================================
\subsection{Experiments and discussion}

Now we are in a position to give some experiments of the different
applications mentionned above. The algorithms induced by our
theoretical analysis involve either solving polynomial systems built from
noisy data or estimating high-dimensional parameters that appear
linearly in equations built also from noisy data. 

The first of these problems is still under very active
research. In our best knowledge, one of the most powerful solver is
FastGb, introduced by Jean-Charles Faugere
\cite{Faugere-98,Faugere-99}. Our experiments involving solving polynomial
systems have been conducted using FastGb. However at this stage, we
will show only synthetic experiments. Future research will be devoted
to the case of real data. The second question is a typical case of
heteroscedastic estimation \cite{Matei-Meer-00} and will be discussed
below. 

{\bf Recovering Epipolar Geometry from a rational cubic and two conics}

We proceed to a synthetic experiment, where the epipolar geometry is
computed from a rational cubic and two conics. The curves are randomly
chosen, as well as the camera.  

Hence the cubic is defined by the following system:
$$
\begin{array}{l}
226566665956452626ZX-1914854993236086169ZT-791130248041963297YZ- \\
1198609868087508022Z^2 + 893468169675527814XT+285940501848919422T^2- \\
179632615056970090YT+277960038226472656Y^2 = 0 \\
555920076452945312XY+656494420457765614ZX-1755155973545148735YZ-\\
1749154450800074954Z^2+ 984240461094724954XT-61309565864179510YT-\\
1802588912007356295ZT+291319745776795474T^2 = 0 \\
1111840152905890624X^2-2905335341664005486ZX-793850352563738017YZ+\\
1286890161434843658Z^2+ 1713207647519936006XT-248798847306328202YT-\\
2942349361064284313ZT+398814386951585134T^2 = 0
\end{array}
$$

The first and the second conic are respectively defined by:
$$
\begin{array}{l}
25X+9Y+40Z+61T = 0 \\
40X^2-78XY+62ZX+11XT+88Y^2+YZ+30YT+81Z^2-5ZT-28T^2
\end{array}
$$
and
$$
\begin{array}{l}
4X-11Y+10Z+57T = 0\\
-82X^2-48XY-11ZX+38XT-7Y^2+58YZ-94YT-68Z^2+14ZT-35T^2 = 0
\end{array}
$$

The camera matrices are given by:
$$
\begin{array}{cc}
{\bf M}_1 = \left[\begin{array}{cccc}
		-87 & 79 & 43 & -66 \\
		-53 & -61 & -23 & -37 \\
		31 &  -34 &  -42 & 88
	    \end{array}\right]
&
{\bf M}_2 = \left[\begin{array}{cccc}
		-76 & -65 & 25 & 28 \\
		-61 & -60 & 9 & 29 \\
		-66 & -32 & 78 & 39
	    \end{array}\right]
\end{array}
$$

Then we form the Extended Kruppa's Equations for each curve. From a
computational point of view, it is crucial to enforce the constraint
that each $\la$ is different from zero. Mathematically this means that
the computation is done in the localization with respect to each
$\la$. 

As expected, we get a zero-dimension variety which degree is one. Thus
there is a single solution to the epipolar geometry given by the
following fundamental matrix: 
$$
{\bf F} = \left[\begin{array}{ccc}
		-\frac{511443}{13426} & -\frac{2669337}{13426} & -\frac{998290}{6713} \\
		\frac{84845}{2329} & \frac{23737631}{114121} & \frac{14061396}{114121} \\
		\frac{1691905}{228242} & \frac{3426650}{114121} & \frac{8707255}{228242}
	  \end{array}\right]
$$

{\bf Reconstruction of a spatial quartic in $\PP^3$}

\begin{figure}
\begin{center}
\psfig{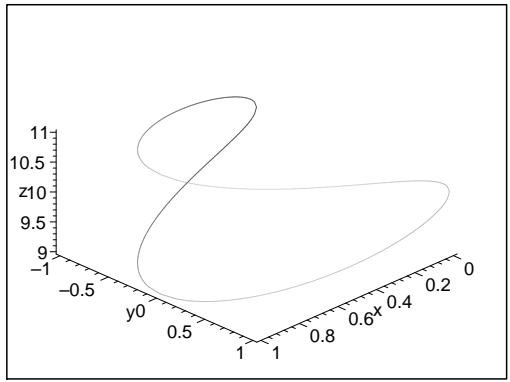}
\end{center}
\caption{A spatial quartic}
\label{spatial_quartic}
\end{figure}

Consider the curve $X$, drawn in figure \ref{spatial_quartic}, defined
by the following equations:  
$$
\begin{gathered} 
F_1(x,y,z,t)=x^2+y^2-t^2 \\ F_2(x,y,z,t)=xt-(z-10t)^2 
\end{gathered} 
$$

The curve $X$ is smooth and irreducible, and has degree $4$ and genus
$1$. We define two camera matrices:      
$$ 
\begin{array}{cc} 
{\bf M}_1=\left[\begin{array}{cccc} 1 & 0 & 0 & 5 \\ 0 & 0 & 1 & -2 \\
0 & -1 & 0 & -10 \end{array}\right] & {\bf  M}_2=\left[\begin{array}{cccc} 1 & 0 & 0 & -10 \\ 0 & 0 & -1 & 0 \\ 0 & 1 & 0 & -10 \end{array}\right]  
\end{array}   
$$

Then the curve is reconstructed from the two projections. As expected
there are two irreducible components. One has degree $4$ and is the
original curve, while the second has degree $12$.

{\bf Reconstruction using the Grassmannian}

\begin{figure}
\begin{center} 
\psfig{figure=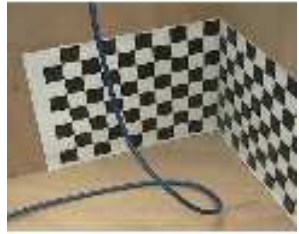,height={5.0cm}} 
\caption{One the six views of an electric thread that were used to perform the reconstruction.}
\label{Thread1}
\end{center}
\end{figure}

\begin{figure}
\begin{center} 
\psfig{figure=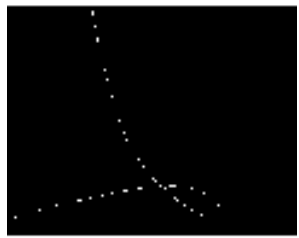,height={5.0cm}} 
\caption{An electric thread after segmentation and thinning.} 
\label{Thread2}
\end{center}
\end{figure}

\begin{figure}[t]
\begin{center} 
\psfig{figure=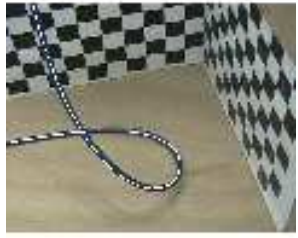,height={5.0cm}}  
\caption{Reprojection on a new image.}
\label{Reproj}
\end{center}
\end{figure}

For the next experiment, we consider six images of an electric wire
--- one of the views  is shown in figure \ref{Thread1} and the image
curve after segmentation and thinning is shown in  figure
\ref{Thread2}. Hence for each of the images, we extracted a set of
points lying on the thread. No fitting is performed in the image
space. For each image, the camera matrix is calculated using the
calibration pattern. Then we proceeded to compute the Chow polynomial
$\Gamma$ of the curve in space. The curve $X$ has degree   
$3$. Once $\Gamma$ is computed, a reprojection is easily  
performed, as shown in figure \ref{Reproj}.

The computation of the Chow Polynomial involves an estimation
problem. Moreover as mentioned above, the Chow polynomial is not
uniquely defined. In order to get a unique solution, we have to add
some constraints to the estimation problem which do not distort the
geometric meaning of the Chow polynomial. This is simply done by
imposing to the Chow polynomial to vanish over $W_d$ additional
arbitrary points of $\PP^5$ which do not lie on $\G(1,3)$. The number
of additional points necessary to get a unique solution is $W_d =
\binom{d+5}{d} - N_d$, where $d$ is the degree of the Chow polynomial.

As we shall the estimation of the Chow polynomial is a typical case of
heteroscedastic estimation. Every 2D measurement ${\bf p}$ is corrupted by
additive noise, which we consider as an isotropic Gaussian
noise ${\cal N}(0,\sigma)$. The variance is estimated to be about 2 pixels. 

For each 2D point ${\bf p}$, we form the optical ray it generates
${\bf L} = \what{\bf M} {\bf p}$. Then the estimation of the
Chow polynomial is made using the optical rays ${\bf L}$. In order
to avoid the problem of scale, the Pl\"ucker coordinates of each line
are normalized such that the last coordinate is equal to one. Hence
the lines are represented by vectors in a five-dimensional affine
space, denoted by ${\bf L}_a$. Hence if $\theta$ is a vector
containing the coefficient of the Chow polynomial $\Gamma$, $\theta$ is
solution of the following problem:
$$
Z({\bf L}_a)^T \theta = 0, \mbox{ for all optical rays},
$$
with $\parallel \theta \parallel = 1$ and $Z({\bf L}_a)$ is a
vector which coordinates are monomials generated by the
coordinates of ${\bf L}_a$. Following
\cite{Chojnacki-all-00,Matei-Meer-00}, in order to obtain a reliable 
estimate, the solution $\theta$ is computed using a maximum likelihood
estimator. This allows to take into account the fact that each $Z({\bf
L}_a)$ has a different covariance matrix, or in other terms that
the noise is {\it heteroscedastic}. More precisely, each $Z({\bf L}_a)$ has the
following covariance matrix:
$$
{\bf C}_{{\bf L}} = {\bf J}_\phi {\bf J}_n \what{\bf M} 
			\left[\begin{array}{ccc}  
				\sigma & 0 & 0 \\
				0        & \sigma & 0 \\
				0        & 0 & 0
			\end{array}\right]
			{\what{\bf M}}^T {\bf J}_n^T {\bf J}_\phi^T,
$$
where ${\bf M}$ is the camera matrix and ${\bf J}_n$ and ${\bf
J}_\phi$ are respectively the Jacobian matrices of the normalization of ${\bf
L}$ and of the map sending ${\bf L}_a$ to $Z({\bf
L}_a)$. That is for ${\bf L}(t)=[L_1,L_2,L_3,L_4,L_5,L_6]^T$, we
have: 
$$
{\bf J}_n = \left[\begin{array}{cccccc}
	\frac{1}{L_6} & 0 & 0 & 0 & 0 & -\frac{L_1}{L_6^2} \\
	0 & \frac{1}{L_6} & 0 & 0 & 0 & -\frac{L_2}{L_6^2} \\
	0 & 0 & \frac{1}{L_6} & 0 & 0 & -\frac{L_3}{L_6^2} \\
	0 & 0 & 0 & \frac{1}{L_6} & 0 & -\frac{L_4}{L_6^2} \\
	0 & 0 & 0 & 0& \frac{1}{L_6}  & -\frac{L_5}{L_6^2} \\
	\end{array}\right],
$$
and ${\bf J}_\phi$ is similarly computed. Then we use the method
presented in \cite{Chojnacki-all-00} to perform 
the estimation. It is worth to note that the estimation is reliable
because the initial guess of the algorithm was well chosen and because
the number of measurements is very large. It is necessary to
use a very large number of measurements for two reasons. First the
dimension of the parameter space is quite high and secondly the measurments are
concentrated on a part of the space (over the Grassmannian $\G(1,3)$).

{\bf Synthetic Trajectory Triangulation}

Let ${\bf P} \in \PP^3$ be a point moving on a cubic, as follows:
$$
{\bf P}(t) = \left[\begin{array}{c}
		t^3 \\ 2t^3+3t^2 \\  t^3+t^2+t+1 \\ t^3+t^2+t+2
	     \end{array}\right]
$$

It is  viewed by a moving camera. At each time instant a picture is
made, we get a 2D point ${\bf p}(t) = [x(t), y(t)]^T = [\frac{{\bf
m}_1^T(t) {\bf P}(t)}{{\bf m}_3^T(t) {\bf P}(t)}, \frac{{\bf
m}_2^T(t) {\bf P}(t)}{{\bf m}_3^T(t) {\bf P}(t)}]^T$, where ${\bf M}^T(t)
= [{\bf m}_1(t), {\bf m}_2(t), {\bf m}_3(t)]$ is the transpose of the camera matrix at
time $t$. 

Then we build the set of optical rays generated by the sequence. The
Chow polynomial is then computed and given below:
$$
\begin{array}{l}
\Gamma(L_1,...,L_6) = -72L_2^2L_3+L_1^3-5L_1L_4L_5- \\ 
18L_1L_3L_6+57L_2L_3L_5+ 48L_2L_4L_5-43L_1L_2L_4- \\
10L_1L_3L_5+ 21L_1L_5L_6-30L_1L_4L_6-108L_2L_3L_6+ \\
41L_1L_2L_5+69L_1L_2L_6-26L_1L_2L_3-36L_2L_4^2- \\
21L_2L_5^2+ 3L_3L_5^2-9L_3^2L_5-12L_4^2L_5+ 6L_4L_5^2+ \\
4L_1^2L_4+20L_2^3-13L_3^3+8L_4^3-L_5^3+108L_2^2L_6-\\
120L_2^2L_5+27L_3^2L_6-25L_1^2L_6+57L_2L_3^2+\\
84L_2^2L_4+7L_1L_3^2-L_1^2L_5+31L_1L_2^2+ \\
5L_1^2L_3+L_1L_5^2-11L_1^2L_2+7L_1L_4^2 
\end{array}
$$ 

\begin{figure}
\begin{center}
\begin{tabular}{c}
\psfig{figure=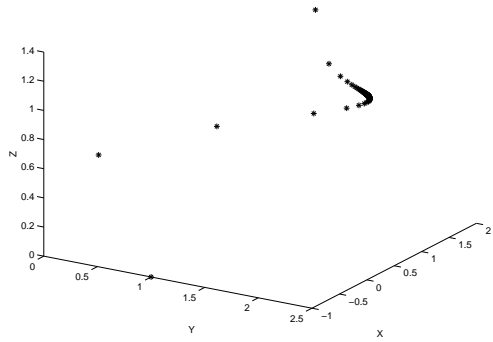,height={5.0cm}}
\end{tabular}
\end{center}
\caption{The 3D locations of the point}
\label{disc}
\end{figure}

From the Chow polynomial, one can extract direclty the locations of
the moving point at each time instant an image was made. This is done
by a two steps computation. The first step consists in giving a
parametric representation of the optical ray generated by the 2D
measurment. During the second step, the pencil of lines passing
through a generic point on the optical ray is considered. For this
generic point to be of the trajectory, the Chow polynomial must vanish
over the pencil. This yields a polynomial system in one variable,
whose root gives the location of the 3D moving point. We show in
figure \ref{disc} the recovered discrete locations of the point in
3D. 

{\bf Trajectory Triangulation from real images}

\begin{figure}
\begin{center}
\begin{tabular}{c}
\psfig{figure=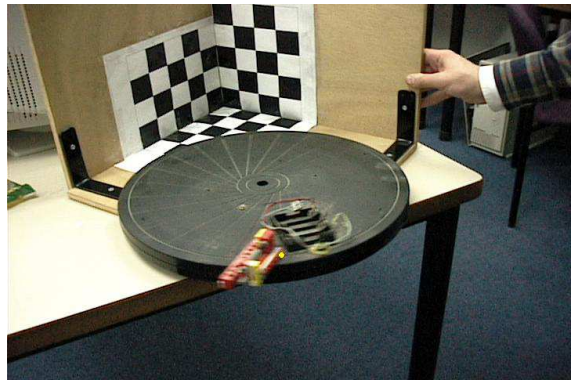,height={5.0cm}}
\end{tabular}
\end{center}
\caption{A moving point over a conic section}
\label{conicIm}
\end{figure}

A point is moving over a conic section. Four static non-synchronized
cameras are looking at it. We show on figure~\ref{conicIm} one image
of one sequence. 

The camera matrices are computed using the
calibration pattern. Every 2D measurement ${\bf p}(t)$ is corrupted by
additive noise, which we consider as an isotropic Gaussian
noise ${\cal N}(0,\sigma)$. The variance is estimated to be about 2
pixels.  

\begin{figure}
\begin{center}
\begin{tabular}{c}
\psfig{figure=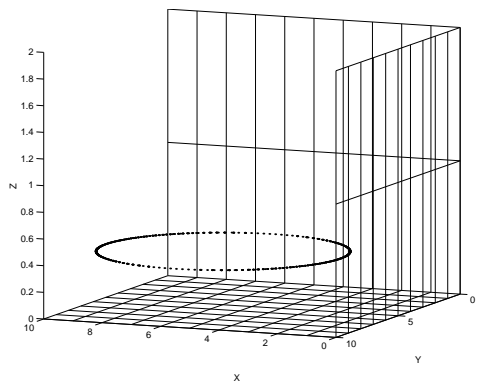,height={5.0cm}}
\end{tabular}
\end{center}
\caption{The trajectory rendered in the calibration pattern}
\label{conicSolution}
\end{figure}

As before the estimation is done from the set of optical rays
generated by the 2D moving point. The estimation is also a case of
heteroscedastic estimation, which was handled with the method
presented in \cite{Chojnacki-all-00}.  The result is stable where
starting with a good initial guess. In order to handle more general situation we further
stabilize it by incorporating some extra constraints that come from
our {\it a-priori} knowledge of the form of the solution. The final
result is presented in figure \ref{conicSolution}.

%========================
%========================


\begin{thebibliography}{}

\bibitem{Avidan-Shashua-00} {\sc S. Avidan and A. Shashua}, 
{\em Trajectory triangulation: 3D reconstruction of moving points from a monocular image sequence},
IEEE Transactions on Pattern Analysis and Machine Intelligence, 22(4):348-357, 2000. 

\bibitem{Barnabei-Brini-Rota-84} {\sc M.Barnabei, A.Brini and G-C.Rota},
{\em On the Exterior Calculus of Invariant Theory},
In Journal of Algebra, 96, 120-160(1985).
 
\bibitem{Berthilsson-all-02} {\sc R. Berthilsson, K. Astrom and
A. Heyden}
{\em Reconstruction of general curves, using factorization and bundle
adjustement},
In International Journal of Computer Vision, 2002.

\bibitem{Chojnacki-all-00} {\sc W. Chojnacki, M. Brooks, A. van den Hengel and D. Gawley},
{\em On the Fitting of Surfaces to Data with Covariances},
PAMI, vol. 22, Nov. 2000.


\bibitem{Faugeras-93} {\sc O.D.~Faugeras} 
{\em Three-Dimensional Computer Vision, A geometric approach}, 
MIT Press, 1993.

\bibitem{Faugeras-Luong-01} {\sc O.D.~Faugeras and Q.T.~Luong}, 
{\em The Geometry of Multiple Images}, MIT Press, 2001.

\bibitem{Faugere-98} {\sc J.C.~Faugere},  
{\em Computing Grobner basis without reduction to zero ($F_5$)},
Technical report, LIP6, 1998.

\bibitem{Faugere-99} {\sc J.C.~Faugere}, 
{\em A new efficient algorithm for computing Grobner basis ($F_4$)}. 

\bibitem{Forsyth} {\sc D. Forsyth}, 
{\em Recognizing algebraic surfaces from their outlines}.

\bibitem{Cross-Zisserman-98} {\sc Cross and A.Zisserman},
{\em Quadric Reconstruction from Dual-Space Geometry}, 1998.

\bibitem{Harris-92} {\sc J.~Harris}, 
{\em Algebraic Geometry, a first course}, 
Springer-Verlag, 1992.

\bibitem{Hartley-Zisserman-00} {\sc R.Hartley and A.Zisserman}, 
{\em Multiple View Geometry in computer vision}, 
Cambridge Univeristy Press, 2000.

\bibitem{Hartshorne-77} {\sc R. Hartshorne}, 
{\em Algebraic Geometry},
Springer, 1977.

\bibitem{Kahl-Heyden-98} {\sc F.~Kahl and A.~Heyden}, 
{\em Using Conic Correspondence in Two Images to Estimate the Epipolar Geometry}, 
In Proceedings of the International Conference on Computer Vision, 1998. 

\bibitem{Kaminski-Shashua-00} {\sc J.Y.~Kaminski and A.~Shashua}, 
{\em On Calibration and Reconstruction from Planar Curves}, 
In Proceedings European Conference on Computer Vision, 2000. 

\bibitem{Ma-Chen-94} {\sc S.D.~Ma and X.~Chen}, 
{\em Quadric Reconstruction from its Occluding Contours}, 
In Proceedings International Conference of Pattern Recognition, 1994.

\bibitem{Matei-Meer-00} {\sc B. Matei and P. Meer}, 
{\em A General Method for Errors-in-variables Problems in Computer
Vision}, 
In {\em Proceedings of IEEE Conference on Computer Vision and Pattern
recognition}, 2000.  

\bibitem{Miranda-91} {\sc Rick Miranda} 
{\em Algebraic Curves and Riemann Surfaces},
American Mathematical Society, 1991.

\bibitem{Quan-96} {\sc L.~Quan},
{\em Conic Reconstruction and Correspondence from Two Views},
In IEEE Transactions on Pattern Analysis and Machine Intelligence, 18(2), February 1996.


\bibitem{Schmid-Zisserman-98} {\sc C.~Schmid and A.~Zisserman},
{\em The Geometry and Matching of Curves in Multiple Views},
In Proceedings European Conference on Computer Vision, 1998.






\end{thebibliography}
\end{document}